\documentclass[12pt]{article}
\usepackage{amsmath, amssymb,latexsym, epsfig, color, pb-diagram}

\newtheorem{theorem}{Theorem}[section]

\newtheorem{question}[theorem]{Question}
\newtheorem{corollary}[theorem]{Corollary}
\newtheorem{definition}[theorem]{Definition}

\newtheorem{lemma}[theorem]{Lemma}

\newcommand{\mc}{\mathcal}
\DeclareMathOperator{\lk}{lk}
\DeclareMathOperator{\cost}{cost}

\newcommand{\p}{^{\prime}}

\begin{document}

\title{Lower Bounds for  Buchsbaum* Complexes}
\author{Jonathan Browder and  Steven Klee \\
\small Department of Mathematics, Box 354350\\[-0.8ex]
\small University of Washington, Seattle, WA 98195-4350, USA,\\[-0.8ex]
\small \texttt{[browder, klees]@math.washington.edu} }
\maketitle

\begin{abstract}
The class of $(d-1)$-dimensional Buchsbaum* simplicial complexes is studied.  It is shown that the rank-selected subcomplexes of a (completely) balanced Buchsbaum* simplicial complex are also Buchsbaum*.  Using this result, lower bounds on the $h$-numbers of balanced Buchsbaum* simplicial complexes are established.  In addition, sharp lower bounds on the $h$-numbers of flag $m$-Buchsbaum* simplicial complexes are derived, and the case of equality is treated. 
\end{abstract}


\section{Introduction}

One commonly studied invariant of a finite, $(d-1)$-dimensional simplicial complex $\Delta$ is its $f$-vector $f(\Delta) = (f_{-1},f_0,\ldots,f_{d-1})$, where $f_i$ denotes the number of $i$-dimensional faces in $\Delta$.  It is equivalent, and oftentimes more convenient, to study the $h$-vector $h(\Delta) = (h_0, \ldots, h_d)$ defined by the relation $\sum_{j=0}^dh_j\lambda^{d-j} = \sum_{i=0}^df_{i-1}(\lambda-1)^{d-i}$.  

When studying the $h$-numbers of simplicial complexes, it is natural to study the class of simplicial complexes that are Cohen-Macaulay over a fixed field $\mathbf{k}$.  This includes the classes of $\mathbf{k}$-homology balls and $\mathbf{k}$-homology spheres.  A more specialized class of simplicial complexes is the class of complexes that are doubly Cohen-Macaulay (2-CM) over $\mathbf{k}$.  This class was introduced and studied by Baclawski \cite{Bac82}.  Any complex that is 2-CM over $\mathbf{k}$ is Cohen-Macaulay over $\mathbf{k}$ with non-vanishing top dimensional homology (computed with coefficients in $\mathbf{k}$).  For example, $\mathbf{k}$-homology spheres are 2-CM over $\mathbf{k}$, but $\mathbf{k}$-homology balls are not 2-CM over $\mathbf{k}$.  

The advantage of studying topological manifolds is that they are locally homeomorphic to topological balls.  Using this local property as motivation, it is natural to define the class of Buchsbaum simplicial complexes.  We say that a simplicial complex $\Delta$ is Buchsbuam over $\mathbf{k}$ if the link of each vertex of $\Delta$ is Cohen-Macaulay over $\mathbf{k}$.  Hence $\mathbf{k}$-homology manifolds are Buchsbaum simplicial complexes over $\mathbf{k}$.  If $\Delta$ is a Buchsbaum simplicial complex, it is convenient to study the $h\p$-numbers of $\Delta$, denoted by $h_j\p(\Delta)$, (defined in Section \ref{notation})  which encode both the $h$-numbers of $\Delta$ and the underlying geometric structure of $\Delta$.  

Athanasiadis and Welker \cite{AthW} define the class of Buchsbaum* complexes that specialize Buchsbaum complexes in the same way that 2-CM complexes specialize Cohen-Macaulay complexes.  They show that if $\Delta$ is Buchsbaum* over $\mathbf{k}$, then the link of each vertex of $\Delta$ is 2-CM over $\mathbf{k}$ and that a homology manifold is Buchsbaum* over $\mathbf{k}$ if and only if it is orientable over $\mathbf{k}$.  

Stanley \cite{S79} introduces the class of \textit{balanced} simplicial complexes and shows that the rank selected subcomplexes of a balanced Cohen-Macaulay complex are Cohen-Macaulay.  This result easily generalizes to the classes of 2-CM and Buchsbaum simplicial complexes.  Our first goal in this paper is to show (Theorem \ref{rank-selected}) that the rank-selected subcomplexes of a Buchsbaum* simplicial complex are Buchsbaum*.  This result answers a question posed in \cite{AthW} in the affirmative. 

Barnette's Lower Bound Theorem \cite{B73} says that if $\Delta$ is a $(d-1)$-dimensional homology manifold without boundary and $d \geq 3$, then $h_2(\Delta) \geq h_1(\Delta)$.  This gives a sharp lower bound on all the $f$-numbers of $\Delta$ in terms of $d$ and $n$, the number of vertices in $\Delta$.  Equality in this lower bound is achieved by a \textit{stacked} simplicial $(d-1)$-sphere on $n$ vertices. 

Kalai \cite{k87} used the theory of rigidity frameworks to prove that $h_2 \geq h_1$ for simplicial homology manifolds without boundary and $d \geq 3$.  Nevo \cite{Ne07} extended this result to the class of 2-CM complexes, and Athanasiadis and Welker \cite{AthW} further extend this result to the class of connected Buchsbaum* complexes.  Using these results, together with Theorem \ref{rank-selected}, we extend a result of Goff, Klee, and Novik \cite{GKN}, showing that $2h_2 \geq (d-1)h_1$  for all balanced, connected Buchsbaum* simplicial complexes with $d \geq 3$.  As in the case of Barnette's LBT, this bound is sharp, with equality achieved by a so-called \textit{stacked cross-polytopal sphere}.  

Athanasiadis and Welker \cite{AthW} show that $h\p_j(\Delta) \geq {d \choose j}$ when $\Delta$ is a flag, Buchsbaum* simplicial complex of dimension $d-1$.  Equality is achieved in this bound by $\mathcal{P}^{\times}_d$, the boundary complex of a $d$-dimensional cross polytope.  We generalize this result, showing that $h\p_j(\Delta) \geq {d \choose j}m^j$ when $\Delta$ is a flag, $m$-Buchsbaum* simplicial complex of dimension $d-1$ in Theorem \ref{m-buchs}.  Moreover, we show that equality holds here only for the $d$-fold join of $m+1$ disjoint vertices, answering a question posed in \cite{AthW} in the case that $m=1$.


The paper is structured as follows.  In Section 2, we review the definitions and background information that will be necessary for the remainder of the paper.  In Section 3, we study balanced Buchsbaum* complexes, proving that the rank selected subcomplexes of Buchsbaum* complexes are Buchsbaum* (Theorem \ref{rank-selected}).  In Section 4, we prove lower bounds on balanced Buchsbaum* complexes (Theorem \ref{BBS-LBT}) and flag Buchsbaum* complexes (Theorem \ref{m-buchs}).  


\section{Notation and Definitions} \label{notation}
 We begin by reviewing some basic definitions on simplicial complexes.  A \textit{simplicial complex} $\Delta$ on vertex set $V=V(\Delta)$ is a collection of subsets $\tau \subseteq V,$ called \textit{faces}, that is closed under inclusion.  We say that a simplicial complex $\Delta$ is \textit{pure} if all of its \textit{facets} (maximal faces under inclusion) have the same dimension.  

The \textit{dimension} of a face $\tau \in \Delta$ is $\dim \tau = |\tau|-1$, and the dimension of $\Delta$ is $\dim \Delta = \max\{\dim \tau: \tau \in \Delta\}$.  The $f$-vector of $\Delta$ is the vector $f(\Delta) = (f_{-1},f_0,\ldots,f_{d-1})$ where $\dim \Delta = d-1$ and the \textit{$f$-numbers} $f_i = f_i(\Delta)$ denote the number of $i$-dimensional faces of $\Delta$.  It is oftentimes more convenient to study the $h$-numbers $h_j(\Delta)$ defined by the relation $\sum_{j=0}^dh_j\lambda^{d-j} = \sum_{i=0}^df_{i-1}(\lambda-1)^{d-i}$.  It is easy to see that the $h$-numbers of $\Delta$ can be written as integer linear combinations of its $f$-numbers and that the $f$-numbers of $\Delta$ can be written as \textit{nonnegative} integer linear combinations of its $h$-numbers.  In particular, bounds on the $h$-numbers of $\Delta$ immediately yield bounds on the $f$-numbers of $\Delta$. 

If $\Delta$ is a simplicial complex and $\tau$ is a face of $\Delta$, the \textit{contrastar} of $\tau$ in $\Delta$ is $\cost_{\Delta}(\tau):= \{\sigma \in \Delta: \sigma \nsupseteq \tau\}$, and the \textit{link} of $\tau$ in $\Delta$ is $\lk_{\Delta}(\tau):= \{\sigma \in \Delta: \sigma \cap \tau = \emptyset, \sigma \cup \tau \in \Delta\}$.  If $A \subset V(\Delta)$ is a collection of vertices in $\Delta$, then $\Delta-A$ is the \textit{restriction} of $\Delta$ to $V(\Delta) \setminus A$.  When $A$ consists of a single vertex $v \in V(\Delta)$, we simply write $\Delta-v$ instead of $\Delta - \{v\}$. If $\Gamma$ and $\Delta$ are simplicial complexes on disjoint vertex sets, their simplicial join is the $(\dim \Gamma + \dim \Delta +1)$-dimensional simplicial complex $$\Gamma * \Delta:= \{\sigma \cup \tau: \sigma \in \Gamma, \tau \in \Delta\}.$$

We are particularly interested in studying the class of balanced simplicial complexes, introduced by Stanley \cite{S79}.

\begin{definition}
A $(d-1)$-dimensional simplicial complex $\Delta$ is \textbf{balanced} if there is a coloring $\kappa: V(\Delta) \rightarrow [d]$ with the property that $\kappa(u) \neq \kappa(v)$ for all edges $\{u,v\} \in \Delta$.  We assume that a balanced complex comes equipped with such a coloring $\kappa$.
\end{definition}

The order complex of a graded poset of rank $d$ is one example of a balanced simplicial complex.  If $\Delta$ is a balanced simplicial  complex and $S \subseteq [d]$, it is often important to study the \textit{$S$-rank selected subcomplex} of $\Delta$, which is defined as the collection of faces in $\Delta$ whose vertices are colored by $S$.  Specifically,
\begin{displaymath}
\Delta_S = \{\tau \in \Delta: \kappa(\tau) \subseteq S\}.
\end{displaymath}
In \cite{S79} Stanley showed that
\begin{equation}\label{hnums}
h_i(\Delta) = \sum_{|S|=i}h_i(\Delta_S);
\end{equation}
and that if $\Delta$ is Cohen-Macaulay, then so are its rank-selected subcomplexes.

A more specialized class of Cohen-Macaulay complexes is the class of \textit{doubly Cohen-Macaulay} (2-CM) complexes.  A $(d-1)$-dimensional simplicial complex $\Delta$ is 2-CM (over $\mathbf{k}$) if it is Cohen-Macaulay and $\Delta-v$ is Cohen-Macaulay of dimension $d-1$ for all vertices $v \in \Delta$.  Walker \cite{W81} showed that double Cohen-Macaulayness is a topological property, meaning that it only depends on the homeomorphism type of the \textit{geometric realization} $|\Delta|$ of $\Delta$. In particular, if a $(d-1)$-dimensional simplicial complex $\Delta$ is 2-CM, then $\cost_{\Delta}\tau$ is Cohen-Macaulay of dimension $d-1$ for all nonempty faces $\tau \in \Delta$. 

Athanasiadis and Welker \cite{AthW} define Buchsbaum* complexes as specializations of Buchsbaum complexes, much in the same sense that 2-CM complexes are specializations of Cohen-Macaulay complexes.  For the purposes of this paper, we will use the following definition of a Buchsbaum* complex.
\begin{definition}
A $(d-1)$-dimensional simplicial complex $\Delta$ that is Buchsbaum over $\mathbf{k}$ is {\bf Buchsbaum*} over $\mathbf{k}$ if, for all $p \in |\Delta|$, the canonical map $\rho_*: \widetilde{H}_{d-1}(|\Delta|;\mathbf{k}) \rightarrow \widetilde{H}_{d-1}(|\Delta|,|\Delta|-p;\mathbf{k})$ is surjective.  
\end{definition}

Henceforth we will fix a field $\mathbf{k}$.  When we say that a simplicial complex $\Delta$ is Buchsbaum* without qualification, we implicitly mean that $\Delta$ is Buchsbaum* over $\mathbf{k}$.  Moreover, we will implicitly compute all homology groups with coefficients in $\mathbf{k}$, and we will suppress this from our notation for convenience.  

First we will give an equivalent definition of the Buchsbaum* property in a combinatorial language.  

\begin{lemma}
Let $\Delta$ be a $(d-1)$-dimensional Buchsbaum simplicial complex. The following are equivalent. 
\begin{enumerate}
\item[\rm{(a)}] $\Delta$ is Buchsbaum*;
\item[\rm{(b)}] For all faces $\sigma \subseteq \tau$ of $\Delta$, the map $$j_*: \widetilde{H}_{d-1}(\Delta,\cost_{\Delta}(\sigma)) \rightarrow \widetilde{H}_{d-1}(\Delta,\cost_{\Delta}(\tau)),$$ induced by inclusion, is surjective;
\item[\rm{(c)}] For all faces $\tau \in \Delta$, the map $$\rho_*:\widetilde{H}_{d-1}(\Delta) \rightarrow \widetilde{H}_{d-1}(\Delta,\cost_{\Delta}(\tau)),$$ induced by inclusion, is surjective.
\end{enumerate}
\end{lemma}
\begin{proof}
The equivalence of (a) and (b) is Proposition 2.8 in \cite{AthW}.  Taking $\sigma = \emptyset$ shows that (b) implies (c).  Next, consider any point $p \in |\Delta|$, and let $\tau$ be the unique minimal face of $|\Delta|$ whose relative interior contains $p$.  Then $|\Delta|-p$ retracts onto $\cost_{\Delta}(\tau)$ and (c) implies (a).
\end{proof}

The $h\p$-vector of a Buchsbaum complex $\Delta$ encodes both the underlying geometry of $\Delta$ and the combinatorial data of $\Delta$.  Let $\widetilde{\beta}_i(\Delta):= \dim_{\mathbf{k}}\widetilde{H}_i(\Delta;\mathbf{k})$ denote the (reduced) \textit{$\mathbf{k}$-Betti numbers} of $\Delta$.  The $h\p$-numbers of $\Delta$ are defined by 
\begin{displaymath}
h_j\p(\Delta) := h_j(\Delta) + {d \choose j} \sum_{i=0}^{j-1}(-1)^{j-i-1}\widetilde{\beta}_{i-1}(\Delta).
\end{displaymath}
We encode the $h\p$-numbers of $\Delta$ into the $h\p$-polynomial $h\p_{\Delta}(t):= \sum_{j=0}^dh\p_{j}(\Delta)t^j.$

Athanasiadis and Welker \cite{AthW} prove a number of very nice properties about Buchsbaum* complexes, which we summarize here. 

\begin{theorem}
Let $\Delta$ be a $(d-1)$-dimensional Buchsbaum* complex with $d \geq 2$.  Then
\begin{enumerate}
\item $\widetilde{H}_{d-1}(\Delta) \neq 0$;
\item $\lk_{\Delta}v$ is 2-CM for all vertices $v \in \Delta$; and
\item $\Delta$ is doubly-Buchsbaum, meaning that $\Delta-v$ is a Buchsbaum complex of dimension $d-1$, for all vertices $v \in \Delta$.
\end{enumerate}
\end{theorem}


\section{Balanced Buchsbaum* Complexes}
It is well known (see, for example, \cite{S79}) that if $\Delta$ is a  balanced, Cohen-Macaulay complex, then $\Delta_S$ is Cohen-Macaulay for any $S \subseteq [d]$.  It is easy to see that this result generalizes to the classes of $2$-CM complexes and Buchsbaum complexes.  The purpose of this section is to generalize this result to the class of Buchsbaum* complexes. 

\begin{theorem} \label{rank-selected}
Let $\Delta$ be a $(d-1)$-dimensional balanced, doubly-Buchsbaum complex.  For any $S \subset [d]$, the rank selected subcomplex $\Delta_S$ is Buchsbaum*.  
\end{theorem}

In particular, since a Buchsbaum* complex is doubly-Buchsbaum, Theorem \ref{rank-selected} implies that the rank selected subcomplexes of a Buchsbaum* complex are Buchsbaum*.

We begin with a series of lemmas, the first of which is well-known (see, for example, \cite{Munkres84}).

\begin{lemma} \label{lem2}
Let $\Gamma$ be a simplicial complex, and let $\tau$ be a nonempty face in $\Gamma$.  Then $\widetilde{H}_i(\Gamma, \cost_{\Gamma}(\tau)) \cong \widetilde{H}_{i-|\tau|}(\lk_{\Gamma}\tau)$.
\end{lemma}

\begin{lemma} \label{lem3}
Let $\Delta$ be a balanced, doubly Buchsbaum simplicial complex of dimension $d-1$, and let $c \in [d]$. Choose vertices $v_1, \ldots, v_i$ of color $c$, and let $\Delta_{i-1} = \Delta \setminus\{v_1, \ldots, v_{i-1}\}$ and $\Delta_i = \Delta \setminus\{v_1, \ldots, v_i\}$.  Then for $S = [d]-c$ and any nonempty face $\tau \in \Delta_S$,
\begin{displaymath}
\widetilde{H}_{d-2}(\cost_{\Delta_{i-1}}(\tau), \cost_{\Delta_i}(\tau)) = 0.
\end{displaymath}
\end{lemma}

\begin{proof}
By Lemma  \ref{lem2},
\begin{eqnarray*}
\widetilde{H}_{d-2}(\cost_{\Delta_{i-1}}(\tau),\cost_{\Delta_i}(\tau)) \cong \widetilde{H}_{d-3}(\lk_{\cost_{\Delta_{i-1}}(\tau)}v_i).
\end{eqnarray*}

Since $\Delta$ is balanced, $\lk_{\cost_{\Delta_{i-1}}(\tau)}v_i = \lk_{\cost_{\Delta}(\tau)}v_i$, and
\begin{displaymath}
\lk_{\cost_{\Delta}(\tau)}v_i = \{\sigma \in \Delta: \sigma \nsupseteq \tau, v_i \notin \sigma, \sigma \cup v_i \in \Delta\} = \cost_{\lk v_i}(\tau).
\end{displaymath}
\noindent Since $\lk_{\Delta}v_i$ is 2-CM of dimension $d-2$, it follows that $\widetilde{H}_{d-3}(\cost_{\lk v_i}(\tau)) = 0$.

\end{proof}

Now we proceed with the proof of Theorem \ref{rank-selected}.

\begin{proof}(Theorem \ref{rank-selected}) \\
Fix a coloring $\kappa: V(\Delta) \rightarrow [d]$.  We need only consider those $S \subset [d]$ with $|S| = d-1$.  Inductively, this is sufficient as any Buchsbaum* complex is doubly Buchsbaum.   Suppose $S = [d]-\{c\}$ and consider the vertices $\{v_1, \ldots, v_k\} \in \Delta$ with $\kappa(v_i)=c$.  For $1 \leq i \leq k$, let $\Delta_i: = \Delta - \{v_1, \ldots, v_{i}\}$, and when $i=0$, set $\Delta_0:= \Delta$.

Let $\tau$ be a nonempty face in $\Delta_S \subset \Delta_{k-1} \subset \cdots \subset \Delta_1 \subset \Delta$.  We claim that for any $0 \leq i \leq k$, the canonical map $$\rho^i_*:\widetilde{H}_{d-2}(\Delta_i) \rightarrow \widetilde{H}_{d-2}(\Delta_i, \cost_{\Delta_i}(\tau))$$ is a surjection.  We proceed by induction on $i$.

When $i=0$, $\Delta_0=\Delta$ is Buchsbaum so $\widetilde{H}_{d-2}(\Delta,\cost_{\Delta}(\tau))=0$, and the map $\rho^0_*$ is surely surjective.

Suppose now that $i>0$ and consider the long exact sequence for the pair $(\cost_{\Delta_{i-1}}(\tau),\cost_{\Delta_{i}}(\tau))$:
\begin{displaymath}
\rightarrow \widetilde{H}_{d-2}(\cost_{\Delta_{i-1}}(\tau),\cost_{\Delta_{i}}(\tau)) \rightarrow \widetilde{H}_{d-3}(\cost_{\Delta_{i}}(\tau)) \stackrel{\iota_*}{\rightarrow} \widetilde{H}_{d-3}(\cost_{\Delta_{i-1}}(\tau)) \rightarrow.
\end{displaymath}
\noindent By Lemma \ref{lem3}, $H_{d-2}(\cost_{\Delta_{i-1}}(\tau),\cost_{\Delta_{i}}(\tau)) \cong 0$, and hence the map $\iota_*$ is an injection. 

Next, we consider the inclusion map $(\Delta_i,\cost_{\Delta_{i}}(\tau)) \hookrightarrow (\Delta_{i-1},\cost_{\Delta_{i-1}}(\tau))$, which induces the following commutative diagram of long exact sequences.

\[
\begin{diagram}
\node{\widetilde{H}_{d-2}(\Delta_i)} \arrow{r,t}{\rho^i_*} \arrow{s} \node{\widetilde{H}_{d-2}(\Delta_i,\cost_{\Delta_{i}}(\tau))} \arrow{s} \arrow{r,t}{\partial_i} \node{\widetilde{H}_{d-3}(\cost_{\Delta_{i}}(\tau))} \arrow{s,r}{\iota_*}  \\
\node{\widetilde{H}_{d-2}(\Delta_{i-1})} \arrow{r,t}{\rho^{i-1}_*} \node{\widetilde{H}_{d-2}(\Delta_{i-1},\cost_{\Delta_{i-1}}(\tau))} \arrow{r,t}{\partial_{i-1}} \node{\widetilde{H}_{d-3}(\cost_{\Delta_{i-1}}(\tau))}
\end{diagram}
\]

By the inductive hypothesis, the map $\rho^{i-1}_*$ is surjective and so the map $\partial_{i-1}$ is the zero map.  By commutativity of the above diagram, $\iota_* \circ \partial_i$ is the zero map, and since $\iota_*$ is an injection, the map $\partial_i$ is also the zero map.  Thus by exactness, $\rho^i_*: \widetilde{H}_{d-2}(\Delta_i) \rightarrow \widetilde{H}_{d-2}(\Delta_i,\cost_{\Delta_{i}}(\tau))$ is a surjection.

This establishes the claim.  In particular, when $i=k$, we have shown that the canonical map $\rho_*: \widetilde{H}_{d-2}(\Delta_S) \rightarrow \widetilde{H}_{d-2}(\Delta_S,\cost_{\Delta_S}(\tau))$ is surjective, and hence $\Delta_S$ is Buchsbaum*.
\end{proof}

\begin{definition} {\rm{(\cite{AthW}, Definition 5.5)}}
Let $\Delta$ be a $(d-1)$-dimensional simplicial complex, and let $m$ be a nonnegative integer.  We say that $\Delta$ is $m$-Buchsbaum* if $\Delta$ is Buchsbaum and $\Delta-A$ is Buchsbaum* of dimension $d-1$ for any subset $A \subset V(\Delta)$ with $|A| < m$.
\end{definition}

\begin{corollary}
Let $\Delta$ be a $(d-1)$-dimensional balanced simplicial complex that is $m$-Buchsbaum* over $\mathbf{k}$.  For any $S \subseteq [d]$, the rank selected subcomplex $\Delta_S$ is $m$-Buchsbaum* over $\mathbf{k}$.
\end{corollary}

\begin{proof}
By Lemma 5.6 in \cite{AthW}, $\Delta$ is $(m+1)$-Buchsbaum over $\mathbf{k}$.  For any subset $A \subseteq V(\Delta_S)$ with $|A| < m$, the complex $\Delta-A$ is doubly Buchsbaum.  Thus $(\Delta-A)_S = \Delta_S-A$ is Buchsbaum* by Theorem \ref{rank-selected}.
\end{proof}


\section{Lower Bounds}
Fix integers $n$ and $d$ such that $d$ divides $n$.  Let $\mathcal{P}^{\times}_d$ denote the boundary complex of the $d$-dimensional cross polytope.  Following \cite{GKN}, define a \textit{stacked cross-polytopal sphere} $\mc{ST}^{\times}(n,d-1)$ by taking the connected sum of $\frac{n}{d}-1$ copies of $\mathcal{P}^{\times}_d$.  In each connected sum, we identify vertices of the same colors so that $\mc{ST}^{\times}(n,d-1)$ is a balanced $(d-1)$-sphere on $n$ vertices.  

Athanasiadis and Welker (\cite{AthW}, Theorem 4.1) prove that if $\Delta$ is a connected, $(d-1)$-dimensional Buchsbaum* complex with $d \geq 3$, then the graph of $\Delta$ is generically $d$-rigid. This generalizes Nevo's result that $h_2(\Delta) \geq h_1(\Delta)$ when $\Delta$ is 2-CM (\cite{Ne07}, Theorem 1.3).  Using Theorem 4.1 from \cite{AthW} in place of Nevo's result and the conclusion of Theorem \ref{rank-selected}, the techniques used to prove Theorem 5.3 in \cite{GKN} give the following Lower Bound Theorem for balanced Buchsbaum* complexes.

\begin{theorem} \label{BBS-LBT}
Let $\Delta$ be a connected, balanced, Buchsbaum* complex of dimension $d-1$ with $d \geq 3$.  Then $d\cdot h_2(\Delta) \geq {d \choose 2}h_1(\Delta)$.  In particular, if  $d$ divides $n = f_0(\Delta)$, then $f_j(\Delta) \geq f_j(\mathcal{ST}^{\times}(n,d-1)$ for all $j$.
\end{theorem}

Hersh and Novik \cite{HN} define the \textit{short simplicial $h$-numbers} of a simplicial complex $\Delta$ as $\widetilde{h}_j(\Delta):= \sum_{v \in \Delta}h_j(\lk_{\Delta}v)$ for $0 \leq j \leq d-1$. Swartz \cite{Sw04} proves that the short simplicial $h$-numbers satisfy 
\begin{equation} \label{eq2}
\widetilde{h}_{j-1}(\Delta) = j\cdot h_j(\Delta) + (d-j+1)h_{j-1}(\Delta).  
\end{equation}
We use this formula, together with Theorem \ref{BBS-LBT} to prove the following theorem. 

\begin{theorem} \label{h3-h1}
Let $\Delta$ be a balanced Buchsbaum* complex of dimension $d-1$ with $d \geq 4$.  Then $d\cdot h_3(\Delta) \geq {d \choose 3}h_1(\Delta)$.
\end{theorem}
\begin{proof}
The link of each vertex $v \in \Delta$ is 2-CM, and hence by Theorem 5.3 in \cite{GKN}, $2h_2(\lk_{\Delta}v) \geq (d-2)h_1(\lk_{\Delta}v)$ for all $v \in \Delta$.  Thus
\begin{eqnarray*}
2\cdot(3h_3(\Delta)+(d-2)h_2(\Delta)) &=&  2 \cdot \widetilde{h}_2(\Delta) \geq (d-2)\widetilde{h_1}(\Delta)\\
&=& (d-2)(2h_2(\Delta)+(d-1)h_1(\Delta)),
\end{eqnarray*}
and the desired result follows.
\end{proof}

Bj\"orner and Swartz \cite{Sw06} have conjectured that the inequality $h_{d-1} \geq h_1$ holds for all 2-CM complexes with $d \geq 3$.  When $d=4$, the conclusion of Theorem \ref{h3-h1} continues to hold without the assumption that $\Delta$ is balanced, establishing that $h_3 \geq h_1$ for all $3$-dimensional Buchsbaum* (and hence 2-CM) complexes.

Athanasiadis and Welker \cite{AthW} prove that if a $(d-1)$-dimensional Buchsbaum* simplicial complex $\Delta$ is flag, then $h\p_{\Delta}(t) \geq (1+t)^d$, where the inequality is interpreted coefficient-wise.  Motivated by a question in \cite{Ath}, they pose the following question.

\begin{question} \label{lb}
{\rm{(\cite{AthW}, Question 6.5(i)) Let $\Delta$ be a $(d-1)$-dimensional flag, Buchsbaum* simplicial complex.  If $h\p_j = {d \choose j}$ for some $1 \leq j \leq d-1$, is $\Delta$ necessarily isomorphic to $\mathcal{P}^{\times}_d$?}}
\end{question}

We will generalize the result of Athanasiadis and Welker to the class of $m$-Buchsbaum* simplicial complexes and answer Question \ref{lb} for this class.  For fixed positive integers $m$ and $d$, we define the simplicial complex $\mathcal{P}(m+1,d)$ to be the $d$-fold join of $m+1$ disjoint vertices.  We remark first that $\mathcal{P}(m+1,d)$ is $(m+1)$-CM and hence $m$-Buchsbaum* by Proposition 5.6 in \cite{AthW}, and second that $\mathcal{P}(2,d) = \mathcal{P}^{\times}_d$. 

\begin{theorem} \label{m-buchs}
Let $\Delta$ be a $(d-1)$-dimensional flag simplicial complex that is $m$-Buchsbaum* over the field $\mathbf{k}$.  Then $h\p_{\Delta}(t) \geq (1 + mt)^d$.  Moreover, if $h\p_j(\Delta) = {d \choose j}m^j$ for some $1 \leq j \leq d-1$, then $\Delta$ is isomorphic to $\mathcal{P}(m+1,d)$. 
\end{theorem}

\begin{proof}
We prove the claim by induction on $d$.  When $d=1$, it is clear that $f_0(\Delta) \geq m+1$, so suppose that $d \geq 2$.

Let $F$ be a $(d-2)$-dimensional face of $\Delta$.  Since $\lk_{\Delta}F$ is $m$-Buchsbaum*, there are at least $m+1$ vertices $v_1, \ldots, v_{m+1}$ in $\lk_{\Delta}F$.  Since $\Delta$ is flag, no two of these vertices $v_i$ are connected by an edge in $\Delta$.  In particular, this means that $\lk_{\Delta-\{v_1, \ldots, v_i\}}(v_{i+1}) = \lk_{\Delta}v_{i+1}$.  Following the techniques of Theorem 1.3 in \cite{Ath} or Corollary 3.3 in \cite{AthW}, 
\begin{eqnarray*}
h\p_{\Delta}(t) &=& h\p_{\Delta-v_1}(t) + t\cdot h_{\lk v_1}(t) \\
&=& h\p_{\Delta-\{v_1,v_2\}}(t) + t\cdot h_{\lk v_2}(t) + t\cdot h_{\lk v_1}(t) \\
&=& \cdots \\
&=& h\p_{\Delta-\{v_1,v_2,\ldots, v_m\}}(t) + t\cdot h_{\lk v_m}(t) + \cdots + t\cdot h_{\lk v_1}(t) \\
&\geq& h_{\lk v_{m+1}}(t) + t\cdot h_{\lk v_m}(t) + \cdots + t\cdot h_{\lk v_1}(t) \qquad {\rm(\dagger)}\\
&\geq&(1+mt)^{d-1} + mt(1+mt)^{d-1} \\
&=& (1+mt)^d.
\end{eqnarray*}

To obtain line ${\rm(\dagger)}$, we use the fact that $\Delta$ is $m$-Buchsbaum* and hence $(m+1)$-Buchsbaum.  Thus $\Delta - \{v_1,\ldots,v_m\}$ is a $(d-1)$-dimensional Buchsbaum complex and  $h\p_i(\Delta - \{v_1, \ldots, v_m\}) \geq h_i(\lk_{\Delta}v_{m+1})$ for all $0 \leq i \leq d$. 

Suppose next that $h_j\p(\Delta) = {d \choose j}m^j$ for some $1 \leq j \leq d-1$.  When $d=2$, the only case to consider is $j=1$. It is easy to see that the complete bipartite graph on two disjoint vertex sets of size $m+1$ is the only flag (i.e. triangle-free) $m$-Buchsbaum* graph with exactly $2(m+1)$ vertices.  

Suppose now that $d \geq 3$.  From the argument used above to show that $h_j\p \geq {d \choose j}m^j$, it follows that $h_j\p(\Delta) = {d \choose j}m^j$ if and only if $h_j(\lk_{\Delta}u) = {d-1 \choose j}m^j$ and $h_{j-1}(\lk_{\Delta}u) = {d-1 \choose j-1}m^{j-1}$ for all vertices $u \in \Delta$.  In particular, one of the numbers $j$ and $j-1$ lies in the set $\{1,2,\ldots, d-2\}$, and so $\lk_{\Delta}u$ is isomorphic to $\mathcal{P}(m+1,d-1)$ for all vertices $u \in \Delta$ by our inductive hypothesis.  

Choose a vertex $u_1 \in \Delta$, and let $\Gamma:= \lk_{\Delta}u_1$.  Let $v$ be a vertex of $\Gamma$.  Then $v$ has $(m+1)(d-1)$ neighbors in $\Delta$, and $(m+1)(d-2)$ of these neighbors lie in $\Gamma$.  Let $u_1, u_2, \ldots, u_{m+1}$ be the neighbors of $v$ in $\Delta$ that do not lie in $\Gamma$.  Now consider a vertex $v\p \in \Gamma$ that is adjacent to $v$.  Then $\lk_{\Delta}\{v,v\p\}$ has $(m+1)(d-2)$ vertices, and $(m+1)(d-3)$ of these vertices lie in $\Gamma$.  The remaining $m+1$ vertices of $\lk_{\Delta}\{v,v\p\}$ are adjacent to $v$, and the only such vertices are $\{u_1,\ldots,u_{m+1}\}$.  Thus $u_iv\p$ is an edge for all $i$.  Since $\Gamma$ is connected, it follows that $u_iw$ is an edge for all $w \in \Gamma$.  We claim that $\Delta$ is connected, and since $\Delta$ is flag, it follows that $\Delta$ is isomorphic to $\mathcal{P}(m+1,d)$.  

Finally, we show that $\Delta$ is connected.  When $j=1$, this is obvious as each connected component of $\Delta$ requires $(m+1)\cdot d$ vertices.  When $j \geq 2$, it is relatively easy to see that $h_j\p(\Delta) = \sum h_j\p(\Delta_t)$, where the sum is taken over all connected components $\Delta_t$ of $\Delta$.  Each connected component of $\Delta$ is $m$-Buchsbaum*, and the claim follows.  
\end{proof}

Taking $m=1$ answers Question \ref{lb}.  We note that the assumption that $\Delta$ is flag in Theorem \ref{m-buchs} can easily be replaced by the assumption that $\Delta$ is balanced to yield the same conclusion.  The following corollary is immediate and very interesting, especially when $m$ is large. 

\begin{corollary}
Let $\Delta$ be a $(d-1)$-dimensional flag  (or balanced) simplicial complex that is $m$-Buchsbaum* over $\mathbf{k}$.  Then $(-1)^{d-1}\widetilde{\chi}(\Delta) \geq m^d$.  In particular, if $\Delta$ is Cohen-Macaulay over $\mathbf{k}$, then $\widetilde{\beta}_{d-1}(\Delta) \geq m^d$. 
\end{corollary}


\section{Concluding Remarks}

Recall that a set of vertices $\tau$ in a simplicial complex $\Delta$ is called a \textit{missing face} if $\tau \notin \Delta$ but $\sigma \in \Delta$ for all $\sigma \subsetneq \tau$.  A flag simplicial complex, for example, only has missing faces of size two.  For any simplicial complex $\Delta$, let $\Delta^{*q}$ denote the simplicial join of $q$ disjoint copies of $\Delta$.  Let $\Sigma^j$ denote the $j$-dimensional simplex and $\partial\Sigma^j$ its boundary complex.  

In \cite{Ne08}, Nevo studies the class of $(d-1)$-dimensional simplicial complexes with no missing faces of dimension larger than $i$.  For fixed integers $d$ and $i$, write $d=qi+r$ where $q$ and $r$ are integers with $1 \leq r \leq i$.  Nevo defines a certain $(d-1)$-dimensional sphere $S(i,d-1)$ with no missing faces of dimension larger than $i$ by $$S(i,d-1):= (\partial\Sigma^i)^{*q}*\partial\Sigma^r.$$

Goff, Klee, and Novik (\cite{GKN}, Theorem 3.1(2)) prove that $h_{\Delta}(t) \geq h_{S(i,d-1)}(t)$ for all $(d-1)$-dimensional simplicial complexes $\Delta$ that are 2-CM with no missing faces of dimension larger than $i$.  This result and Theorem \ref{m-buchs} motivate the following question.  

Fix integers $m, i,$ and $d$ and consider the class of $(d-1)$-dimensional simplicial complexes that are $m$-CM with no missing faces of dimension larger than $i$.  As before, write $d = qi+r$ with $1 \leq r \leq i$, and consider the simplicial complex
\begin{displaymath}
S(m,i,d-1):= (\text{Skel}_{i-1}(\Sigma^{m+i-2}))^{*q}*\text{Skel}_{r-1}(\Sigma^{m+r-2}).
\end{displaymath}
Notice that $S(2,i,d-1)$ is Nevo's $S(i,d-1)$, and $S(m,1,d-1)$ is $\mathcal{P}(m,d)$.  Each join-summand $\text{Skel}_{i-1}(\Sigma^{m+i-2})$ is $m$-CM, and hence $S(m,i,d-1)$ is a $(d-1)$-dimensional simplicial complex that is $m$-CM with no missing faces of dimension larger than $i$.  It seems natural, therefore, to pose the following question.

\begin{question}
{\rm{Let $\Delta$ be a $(d-1)$-dimensional simplicial complex that is $m$-CM with no missing faces of dimension larger than $i$.  Is it necessarily true that $$h_{\Delta}(t) \geq h_{S(m,i,d-1)}(t)?$$}}
\end{question}


\section*{Acknowledgements}
We are grateful to Christos Athanasiadis and Volkmar Welker for introducing us to the problems discussed in this paper.  Our thanks also go to Christos Athanasiadis and Isabella Novik for a number of helpful conversations during the development of this paper.  Jonathan Browder's research is partially supported by VIGRE NSF Grant DMS-0354131.

\bibliographystyle{plain}
\bibliography{../biblio}
\end{document}